\documentclass{amsart}
\usepackage{amssymb}
\usepackage[frame,ps,dvips,matrix,arrow,curve]{xy}
\input{xypic}

\newtheorem{defi}[equation]{Definition}
\newtheorem{lemma}[equation]{Lemma}
\newtheorem{prop}[equation]{Proposition}
\newtheorem{theor}[equation]{Theorem}

\numberwithin{equation}{section}
\title{The Hochschild cohomology of a Poincar\'{e} algebra}
\author{Po Hu}

\begin{document}

\begin{abstract}
In this note, 
we define the notion of a cactus set, and show that its geometric realization has 
a natural structure as an algebra over Voronov's cactus operad, which is 
equivalent to the framed $2$-dimensional little disks operad
$\mathcal{D}_2$. Using this, we 
show that for a Poincar\'{e} algebra $A$, its Hochschild cohomology 
is an algebra over the (chain complexes of) $\mathcal{D}_2$.
\end{abstract}

\maketitle

\footnote{This material is based upon work supported by the National Science Foundation under Grant No. 0503814.}

In~\cite{cs}, Chas and Sullivan considered $H_{\ast}(LM)$, the integral
singular homology of the free loop space on a compact smooth oriented manifold $M$,
and showed that it has the structure of a Batalin-Vilkovisky algebra, i.~e.
an algebra over the framed 
2-dimensional little disks operad $\mathcal{D}_2$. In this note, we consider
the question of what algebraic structures have the property that its 
cohomology has the structure of a $\mathcal{D}_2$-algebra. The main 
result of the note is that if $A$ is a Poincar\'{e} algebra, then 
the dual of the 
Hochschild homology of $A$ has the natural structure of an algebra 
over the framed 2-dimensional little disks operad $\mathcal{D}_2$. 
Here, by a Poincar\'{e} algebra, we mean a Frobenius algebra that is 
commutative, i.~e. a commutative algebra $A$ with an augmentation 
map into the base field, such that the adjoint map from $A$ to the 
dual of $A$ is an isomorphism. To prove the theorem, we
make use of Voronov's cactus operad~\cite{voronov}, and define 
a structure of \emph{cactus sets}, then show that the Hochschild cohomology 
of a Poincar\'{e} algebra has such a structure. 
This gives a generalization of the string topology result of Chas and 
Sullivan on $H_{\ast}(LM)$. Although strictly speaking,
theirs is not a Poincar\'{e} algebra, but only one up 
to homotopy.

In the first section, we recall the cactus operad
and define cactus sets, and show that 
the simplicial realization of a cactus set has the structure 
of an algebra over the cactus operad. In Section 2, 
we apply this notion to the case of the Hochschild cohomology of a Poincar\'{e} 
algebra. 

\section{Cactus objects}

In this section, we will define the notion of a cactus object
in a symmetric monoidal category with simplicial realization, which is 
a cyclic object in the sense of Connes~\cite{connes}
with certain extra structures, 
and show that the simplicial realization of a cactus object has a natural 
action by the cactus operad defined by Voronov~\cite{voronov}.

In~\cite{voronov}, Voronov defined the cactus operad $\mathcal{C}= \{ 
\mathcal{C}(n) |\ n \geq 1 \}$ as follows. 
For each $n$, an element 
of $\mathcal{C}(n)$ is an ordered configuration of 
$n$ parametrized circles (the ``lobes'' of the cactus), with varying positive
radii, such that the sum of all the radii is $1$, and that 
the dual graph of the configuration is a tree. There are also the 
following additional data:

1. A cyclic ordering of the lobes at each intersection 
point of the circles. 

2. A chosen distinguished point $0_i$ on the $i$-th circle of the 
cactus, for $i = 1, \ldots, n$.

3. A chosen distinguished point $0$ for the entire configuration. If
$0$ is an intersection point of circles, there is also a choice of 
the particular circle on which it lies.  

For instance, the following configuration is an element of $\mathcal{C}(6)$. The 
solid dots are the distinguished points of the lobes, and the $\times$ mark is the 
distinguished point for the entire cactus.

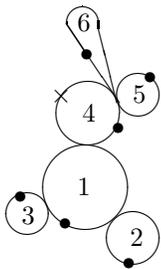
\begin{picture}(150, 130)(-50, 0)
\put(50, 50){\circle{30}}
\put(47, 47){$1$}
\put(68, 31){\circle{20}}
\put(67, 28){$2$}
\put(28, 40){\circle{16}}
\put(26, 36){$3$}
\put(51, 78){\circle{25}}
\put(49, 75){$4$}
\put(70, 85){\circle{17}}
\put(68, 82){$5$}
\put(37, 82){$\times$}
\put(40, 34){$\bullet$}
\put(64, 19){$\bullet$}
\put(23, 44){$\bullet$}
\put(60, 70){$\bullet$}
\put(72, 89){$\bullet$}
\put(63, 82){\line(-2, 3){20}}
\put(62, 82){\line(-1, 4){7}}
\put(49, 112){\oval(12, 13)[t]}
\put(55, 110){\line(0, 1){3}}
\put(47, 109){$6$}
\put(48, 98){$\bullet$}
\end{picture}

The topology on $\mathcal{C}(n)$ is obtained as a quotient of 
$(S^1)^n$. 

The operad structure on $\{ \mathcal{C}(n) \}$
is as follows. Given such a configuration $c \in \mathcal{C}(n)$, 
the choice of $0$ and the 
cyclic ordering of the lobes at each intersection point 
defines a continuous map 
\[ f_c: S^1 \rightarrow c . \] 
Namely, given the standard circle with radius $1$, we can
wrap it around the configuration $c$, starting at the distinguished
point $0$ of $c$, in the manner 
prescribed by the cyclic ordering of the components at each intersection point:
namely, whenever we arrive at an intersection point, we always continue 
onto the circle that comes after the one through which we arrived, in the 
cyclic ordering (note that the circles are parametrized, so 
they come with orientations). Given cacti 
$c$ in $\mathcal{C}(n)$ and $d_i \in \mathcal{C}(k_i)$ for $i =1, 
\ldots, n$, the composition $\gamma(c; d_1, \ldots, d_n) \in 
\mathcal{C}(k_1 + \cdots + k_n)$ 
is obtained by collapsing the $i$-th lobe of $c$ to $d_i$ via this map
(with its chosen distinguished point $0_i$ identified with the distinguished 
point for $d_i$). The distinguished point $0_c$ for the cactus $c$ carries along 
on the lobe to which it belongs. 
\vspace{2mm}

Voronov proved the following theorem. 

\begin{theor}[\cite{voronov}]
The cactus operad is naturally homotopy equivalent to $\mathcal{D}_2$, the 
framed $2$-dimensional little disks operad. 
\end{theor}

\noindent\textbf{Remark:}
Note that as defined, the cactus operad is unbased: there is no $\mathcal{C}(0)$. 
Hence, a $\mathcal{C}$-algebra is a structure that is non-unital. 
However, from the point of view of $\mathcal{D}_2$, we can think of 
$\mathcal{D}_2(0)$ as consisting of a single elemnt, which is a solid little disk 
with no framed little disks inside of it. Composition with this element fills
in one framed little disk. (This is also the point of view from conformal field theory, 
where the $0$-th space of the operad ``caps off'' one inbound boundary component 
of the surface.) From this, one could say that the $\mathcal{C}(0)$ should also 
consist of one configuration, which is a single point. Composition with this element 
contracts one lobe of the cactus to a single point (erasing the distinguished point 
of that lobe). However, we shall not make use of this in this note. 
\vspace{2mm}

Another ingredient which we will need is the notion of 
a cyclic set~\cite{connes}. 
Namely, the cyclic category $\Lambda$ has objects 
all sets $[n] = \{ 0, \ldots, n \}$,
same as the simplicial category $\Delta$, 
and the morphisms are generated by the usual face maps $d_i$ and 
degeneracy maps $s_i$ of the simplicial category with 
the usual relations, as well as one extra 
degeneracy map 
\[ s_{n+1} : [n+1] \rightarrow [n] \]
for each $n \geq 1$, which has the relation 
\[ (d_0 s_{n+1})^{n+1} = Id : [n] \rightarrow [n] . \]
In particular, in the 
cyclic category the number of degeneracy maps $[n+1] \rightarrow 
[n]$ and the number of face maps  
$[n] \rightarrow [n+1]$ are the same.
A cyclic set is a functor 
$S_{\bullet}: \Lambda^{op} \rightarrow Sets$. In other words, 
a cyclic set is a simplicial set $S_{\bullet}$
with an extra degeneracy $S_{n} \rightarrow S_{n+1}$ between the $n$-the 
and the $0$-th simplicial coordinates, for each simplicial degree $n$.
In particular, if $S_{\bullet}$ is 
a cyclic set, then its simplicial realization $|S_{\bullet}|$ has naturally 
an action of $S^1$ (see~\cite{jones}). 

One important property of the cyclic category is that that 
\[ \Lambda \simeq \Lambda^{op} \]
by reversing the faces and degeneracies (see~\cite{connes, elmendorf}). 
(This equivalence is not canonical, 
since one can compose it with any automorphism of $\Lambda$ obtained 
by rotation.) Hence, we also have 
\[ \Lambda^{op}Sets \simeq (\Lambda^{op}Sets)^{op} \]
as categories, and the dual of a cyclic set is again a cyclic set. 

The notion of a cactus object can be thought of as a generalization of a
cyclic object, but with the cactus configurations of $\mathcal{C}$ 
taking the place of the circle. 
For a cactus $c \in \mathcal{C}(n)$, we will need to consider the intersection 
points of circles in $c$. Given an intersection point $x$ of $c$, let the 
\emph{multiplicity} of $x$ denote the number of circles that $x$ lies on.

We will need to consider only the combinatorial 
part of the structure of the cactus operad, and define the following
notion of a ``spiny cactus''. 
For each cactus $c \in \mathcal{C}(n)$, we saw above
that there is a 
well-defined map $f_c: S^1 \rightarrow c$. 
Given $c \in \mathcal{C}(n)$ and positive numbers $j_1, \ldots, j_n$, we
define the set
\[ c_{(j_1-1, \ldots, j_n-1)} \]
to consist of all configurations $(c, X)$, where $X$ is a 
set of chosen points on $c$, such that the points $0$, $0_i$ for $i = 1, 
\ldots, n$, as well as all the intersection points are in $X$,
and there are exact $j_i$ points on the $i$-th circle. 
An intersection point
of multiplicity $n$ is considered to be a point on every circle which
contains the point, hence, it is counted $n$ times. 
Note that the set of points 
on each lobe of the cactus comes with a cyclic ordering, starting with 
$0_i$. With the provision
that an intersection point appears $n$ times, 
there is also a cyclic ordering on the set $X$ of points on the 
entire cactus (having $\Sigma j_i$ points in all), 
arising from the cyclic ordering of the preimages of 
the points on $S^1$ and starting from the distinguished point $0_c$. 
Let $\mathcal{C}(n)_{(j_1, \ldots, j_n)}$ be the disjoint union of 
$c_{(_1, \ldots, j_n)}$ for all $c \in \mathcal{C}(n)$. We make the 
following identifications on $\mathcal{C}(n)_{(j_1, \ldots, j_n)}$. 
For two configurations $(c_1, X_1) \in (c_1)_{(j_1, \ldots, j_n)}$ and 
$(c_2, X_2) \in (c_2)_{(j_1, \ldots, j_n)}$, 
we identify $(c_1, X_1)$ and $(c_2, X_2)$ if $f_{c_2}$ can be 
obtained from $f_{c_1}$ by a 
continuous reparametrization of $S^1$, which takes the points of 
$X_1$ to the points of $X_2$, matching the distinguished points and 
intersection points exactly. This is clearly an equivalence relation, 
and we denote by 
\[ \overline{c}_{(j_1, \ldots, j_n)} \]
the equivalence class of $c_{(j_1, \ldots, j_n)}$. For instance, a representative
of such a ``spiny cactus'' configuration $(c, X)$, where $c \in \mathcal{C}(6)$ 
is the cactus pictured above, looks as follows:

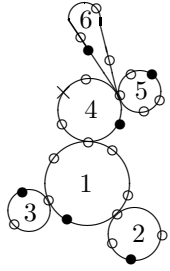
\begin{picture}(150, 130)(-50, 0)
\put(50, 50){\circle{30}}
\put(47, 47){$1$}
\put(68, 31){\circle{20}}
\put(67, 28){$2$}
\put(28, 40){\circle{16}}
\put(26, 36){$3$}
\put(51, 78){\circle{25}}
\put(49, 75){$4$}
\put(70, 85){\circle{17}}
\put(68, 82){$5$}
\put(37, 82){$\times$}
\put(40, 34){$\bullet$}
\put(64, 19){$\bullet$}
\put(23, 44){$\bullet$}
\put(60, 70){$\bullet$}
\put(72, 89){$\bullet$}
\put(69, 75){$\circ$}
\put(64, 89){$\circ$}
\put(75, 79){$\circ$}
\put(60, 81){$\circ$}
\put(59, 36){$\circ$}
\put(76, 28){$\circ$}
\put(33, 41){$\circ$}
\put(20, 32){$\circ$}
\put(57, 25){$\circ$}
\put(48, 63){$\circ$}
\put(35, 57){$\circ$}
\put(57, 60){$\circ$}
\put(63, 43){$\circ$}
\put(38, 70){$\circ$}
\put(45, 103){$\circ$}
\put(51, 114){$\circ$}
\put(56, 94){$\circ$}
\put(47, 87){$\circ$}
\put(63, 82){\line(-2, 3){20}}
\put(62, 82){\line(-1, 4){7}}
\put(49, 112){\oval(12, 13)[t]}
\put(55, 110){\line(0, 1){3}}
\put(47, 109){$6$}
\put(48, 98){$\bullet$}
\end{picture}

Here, the empty circles are the marked points other than the distinguished 
points of the cactus and of the individual lobes. 

In essence, this identification removes the geometric information contained in 
$\mathcal{C}(n)$, and retains only the combinatorial informaion. 
In particular, 
two configurations $(c_1, X_1)$ and $(c_2, X_2)$ are identified
if one can be obtained from the other by changing the radii of the lobes, 
by moving one lobe of the cactus (or any of the 
chosen points) along the circumference of another lobe. However, one 
is not allowed to move a lobe or any of the distinguished points past 
each other, or past any other point of $X$.

For each $1 \leq k \leq n$, and $0 \leq i \leq j_k-1$, 
we have a $(k, i)$-th cyclic degeneracy 
\[ \overline{c}_{(j_1-1, \ldots, j_k-1, \ldots, j_n-1)} \rightarrow 
\overline{c}_{(j_1-1, \ldots, j_k, \ldots, j_n-1)} \]
which is obtained by inserting a new point in the $i$-th position
between two adjacent points on the $k$-th circle. 
The $(k, i)$-th cyclic face map 
\[ \overline{c}_{(j_1-1, \ldots, j_k-1, \ldots, j_n-1)} \rightarrow 
\overline{c^{\prime}}_{(j_1-1, \ldots, j_k -2, \ldots, j_n-1)} \]
is obtained by pinching together two adjacent points at the $i$-th and 
$(i +1)$-st positions on the $k$-th circle of the cactus. Note that 
a cyclic degeneracy always take a spiny cactus based on $c$ to a spiny 
cactus based on the same cactus $c$. However, a cyclic face degeneracy 
map may actually change the cactus configuration itself, if both 
points being pinched together are ``special points'' of the cactus, 
i.~e. distinguished points or intersection points of the lobes. 
For instance, given the first ``spiny'' cactus configuration pictured below,
if we pinch the intersection point of lobes $1$ and $2$ together with the 
intersection point of lobes $1$ and $3$, we get the second configuration,
which is a different cactus. 

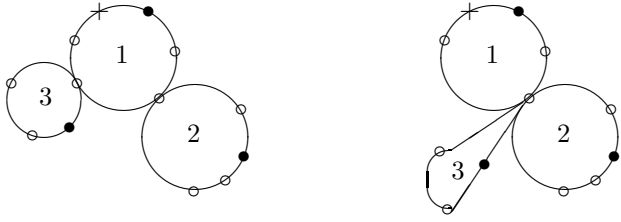
\begin{picture}(150, 130)(-20, 0)
\put(20, 80){\circle{60}}
\put(17, 78){$1$}
\put(47, 50){\circle{40}}
\put(44, 48){$2$}
\put(-10, 64){\circle{27}}
\put(-12, 62){$3$}
\put(27, 95){$\bullet$}
\put(63, 40){$\bullet$}
\put(-3, 51){$\bullet$}
\put(7, 95){$+$}
\put(37, 80){$\circ$}
\put(31, 62){$\circ$}
\put(0, 68){$\circ$}
\put(56, 31){$\circ$}
\put(-1, 84){$\circ$}
\put(-25, 68){$\circ$}
\put(-17, 48){$\circ$}
\put(62, 58){$\circ$}
\put(44, 27){$\circ$}
\put(160, 80){\circle{50}}
\put(157, 78){$1$}
\put(187, 50){\circle{40}}
\put(184, 48){$2$}
\put(167, 95){$\bullet$}
\put(203, 40){$\bullet$}
\put(147, 95){$+$}
\put(177, 80){$\circ$}
\put(139, 84){$\circ$}
\put(202, 58){$\circ$}
\put(184, 27){$\circ$}
\put(196, 31){$\circ$}
\put(171, 62){$\circ$}
\put(171, 62){\line(-2, -3){27}}
\put(171, 63){\line(-3, -2){27}}
\put(144, 34){\oval(18, 22)[l]}
\put(144, 34){$3$}
\put(154, 37){$\bullet$}
\put(140, 20){$\circ$}
\put(137, 42){$\circ$}
\end{picture}

The composition of cacti also translates to this model. Given 
cacti $c \in \mathcal{C}(n)$ and $d_i \in \mathcal{C}(m_i)$ for 
$1 \leq i \leq n$, consider the sets 
$\overline{c}_{(j_1-1, \ldots, j_n-1)}$
and  
\[ \overline{d_1}_{(r_{1,1}-1, \ldots, r_{1,m_1}-1)}, \ldots, 
\overline{d_n}_{(r_{n,1}-1, \ldots, r_{n, m_n}-1)}. \] 
We require that 
for each $1 \leq i \leq n$, 
\[ \Sigma_{l =1}^{m_i} r_{i, m_l} = j_i . \]
Then we have a well-defined composition 
\begin{equation*}
\begin{split} 
\overline{c}_{(j_1-1, \ldots, j_n-1)} \times & 
\overline{d_1}_{(r_{1,1}-1, \ldots, r_{1,m_1}-1)}
\times \cdots \times \overline{d_n}_{(r_{n,1}-1, \ldots, r_{n, m_n}-1)} \\
& \rightarrow 
\overline{\gamma(c; d_1, \ldots, d_n)}_{(s_1-1, \ldots, s_{\Sigma m_i}-1)} 
\end{split}
\end{equation*}
where $\gamma(c; d_1, \ldots, d_n)$ is the composition of the cacti in 
$\mathcal{C}$, and $s_1, \ldots, s_{\Sigma m_i}$ is a permutation of 
$r_{1, 1}, \ldots, r_{n, m_n}$ obtained by the ordering of lobes on the new
cactus. 
In this sense, $\overline{c}_{(\ast, \ldots, \ast)}$ give a partial ``cyclic 
model'' of the cactus $c$. 

Let $(\mathcal{S}, \otimes)$ be a symmetric monoidal category, which also 
has simplicial realization from $\Delta^{op} \mathcal{S}$ to $\mathcal{S}$.
We have the following definition of a cactus object in $\mathcal{S}$. 

\begin{defi}
A cactus object $S_{\bullet}$ in $\mathcal{S}$ is a cyclic 
object in $\mathcal{S}$, with additionally a structure 
map for each $(c, X) \in \overline{c}_{j_1, \ldots, j_n}$ 
\[ \mu_{c, X}: S_{j_1-1} \otimes S_{j_2-1} \otimes \cdots \otimes S_{j_n-1} \rightarrow 
S_{j_1 + \cdots + j_n - 1} . \]
These maps are compatible with the face and degeneracy maps of the cyclic 
set in the manner prescribed by the cactus $c$. 
Namely, there is a cyclic ordering of 
the segments of $X$ over the entire cactus, and the 
total simplicial degree is exactly $j_1 + \cdots j_n -1$.
For $1 \leq k \leq n$ and $0 \leq i \leq j_k-1$, let $f(k, i)$ be the position 
in this overall cyclic ordering corresponding to the $i$-th segment of the 
$k$-th lobe. For the $(k,i)$-th 
degeneracy, we require the following diagram to commute: 
\[ \diagram 
S_{j_1-1} \otimes \cdots \otimes S_{j_k-1} \otimes \cdots \otimes 
S_{j_n-1} \dto_{\sigma_i} \rto^{\mu_{c, X}} & 
S_{j_1 + \cdots + j_n - 1} \dto^{\sigma_{f(k, i)}} \\
S_{j_1-1} \otimes \cdots \otimes S_{j_k} \otimes \cdots \otimes 
S_{j_n-1} \rto_{\mu_{c, X^{\prime}}} & S_{j_1 + \cdots + j_n - 1} 
\enddiagram \]
where $\sigma_i$ is the $i$-th degeneracy. Similarly, for the $(k, i)$-th
face map for $(c, X)$, the following diagram commutes: 
\[ \diagram 
S_{j_1-1} \otimes \cdots \times S_{j_k-1} \otimes \cdots \otimes 
S_{j_n-1} \dto_{\delta_i} \rto^{\mu_{c, X}} & 
S_{j_1 + \cdots + j_n - 1} \dto^{\delta_{f(k, i)}} \\
S_{j_1-1} \otimes \cdots \otimes S_{j_k-2} \otimes \cdots \otimes 
S_{j_n-1} \rto_{\mu_{c^{\prime}, X^{\prime}}} & S_{j_1 + \cdots + j_n - 1} 
\enddiagram \]
where $\delta_i$ denotes the $i$-th face map. 

Further, the maps $\mu_{c, X}$ are compatible with the 
composition of cacti in the following manner. Given $(c, X)$ 
as above, for each $k$ between $1$ and $n$, suppose we also have
$(d_k, Y_k) \in \overline{d_k}_{(r_{k,1}-1, \ldots, r_{k, m_k}-1)}$,
such that $\Sigma r_{k, i} = j_k$. Then we require the following 
diagram to commute: 
\[ \diagram 
(S_{r_{1,1}-1} \otimes \cdots \otimes S_{r_{1, m_1}-1}) \otimes \cdots 
\otimes (S_{r_{n,1}-1} \otimes \cdots \otimes S_{r_{n, m_n}-1}) \drto \dto & \\
S_{j_1-1} \otimes \cdots \otimes S_{j_n-1} \rto & S_{\Sigma j_k -1} . 
\enddiagram \] 
\end{defi}

For simplicity, we state and prove
the following propostion in the context of sets. For 
the category of chain complexes, which will be the case relevant in the
next section, the argument of the proof goes through by taking the chain 
complexes of the standard simplices (and the cactus operad). 

\begin{prop}
If $S_{\bullet}$ is a cactus set, then the simplicial realization of 
$S_{\bullet}$ has naturally an action of the cactus operad $\mathcal{C}$. 
\end{prop}

\begin{proof}
For a cactus $c \in \mathcal{C}(n)$, 
we will define a map 
\[ m_{c}: \Delta^{j_1-1} \times \cdots \times \Delta^{j_n-1} \rightarrow 
\Delta^{j_1 + \cdots j_n -1 +m}  \] 
where $m$ is 
$1 +$ the total of the multiplicities of the intersection points of $c$. 
Namely, for each $1 \leq k \leq n$, let $r_k$ be the radius of 
the $k$-th lobe of the cactus. Given a point with barycentric coordinates 
$(s_{k, 1}, \ldots, s_{k, j_k}) \in \Delta^{j_k -1}$, with 
$s_{k, 1} + \cdots + s_{k, j_k} = 1$, by scaling $r_k$ so that 
the circumference of the $k$-th lobe is $1$,  
starting from the distinguished point of the $k$-th lobe, $(s_{k, 1}, 
\ldots, s_{k, j_k})$ determine $j_k$ points on the $k$-th lobe, which 
divide the circle into $j_k$ subintervals having lengths 
$s_{k,1}, \ldots, s_{k, j_k}$. For each circle, we need to make sure
that any intersection points are among the dividing points; we ensure this 
by using degeneracy maps $\Delta^{j_k} \rightarrow \Delta^{j_k +1}$ to 
add the intersection points on each circle. (Note that hence, each intersection 
point of $c$ is added as many times as its multiplicity. If on a particular 
circle, it happens to be a dividing point from the 
original barycentric coordinates, 
it is still added, and we will get a point in a higher simplex some of whose 
barycentric coordinates are $0$.) Similarly, we also add the distinguished 
point $0$ of the cactus, on the lobe on which it lies. 
Starting at the distinguished point $0$ for the entire
cactus, and using the well-defined flow (which has degree 1 on each circle)
on the cactus, we get a well-defined subdivision of a single circle 
with circumference 1 into $j_1 + \cdots j_k -1$ subintervals.  
The length of these, in cyclic order, are 
$t_1, \ldots, t_{j_1 + \cdots + j_n}$, with $t_1 + \cdots + t_{j_1 + \cdots 
+ j_n} = 1$. Hence,
\[ (t_1, \ldots, t_{j_1 + \cdots +j_n}) \]
is a point in $\Delta^{j_1 + \cdots + j_n -1}$.
Note that this map uses the actual cactus $c \in \mathcal{C}(n)$, including 
the radii of its circles and the exact locations of its distinguished points, 
instead of just its class in $\overline{\mathcal{C}}(n)$.
It is straightforward to see that this map is associative, from the associativity 
of the composition of cacti in $\mathcal{C}$, and is compatible with the 
cyclic set structure of rotation. 

Now for a cactus set $S_{\bullet}$, for each configuration $c \in 
\mathcal{C}(n)$, and $j_1, \ldots, j_n$, we have 
\[ S_{j_1 -1} \times \cdots \times S_{j_n -1} \rightarrow S_{j_1 + \cdots 
+ j_n +m -1} \]
and 
\[ \Delta^{j_1 -1} \times \cdots \times \Delta^{j_n -1} \rightarrow 
\Delta^{j_1 + \cdots + j_n +m -1} . \]
Here, $m$ is $1+$ the total of the multiplicities of the intersection points of $c$,
and the map on $S_{\bullet}$ is obtained from the structure map of a cactus set
by inserting degeneracies at the appropriate positions. 

Now recall that the simplicial realization of $S_{\bullet}$ is given by 
the coequalizer 
\[ \amalg_{m \rightarrow l} S_{l} \times \Delta^m \stackrel{\rightarrow}{\rightarrow}
\amalg_{l} S_l \times \Delta^n \rightarrow |S_{\bullet}| . \]
For $c \in \mathcal{C}(n)$, we have a map
\[ (\amalg_{l} S_l \times \Delta^l)^n \rightarrow \amalg_l S_l \times 
\Delta^l . \]
Namely, consider the source as the sum over all $l_1, \ldots, l_n$ of 
\[ (S_{l_1-1} \times \cdots \times S_{l_n -1}) \times 
(\Delta^{l_1 -1} \times \cdots \times \Delta^{l_n -1}) . \]
As above, 
an element of $\Delta^{l_1 -1} \times \cdots \times \Delta^{l_n -1}$ determines
a set $X$ of distinguished points on $c$. We have the map $m_c$, which 
takes it to $\Delta^{l_1 + \cdots + l_n -1 +m}$. Further, using the 
spiny cactus $(c, X)$, we get a map from 
\[ S_{l_1 -1} \times \cdots \times S_{l_n-1} \rightarrow 
S_{l_1 + \cdots + l_n-1 +m} \]
(after inserting degeneracies at the appropriate places of each 
$S_{l_i -1}$, corresponding to the positions of the intersection points 
on each lobe of the cactus and to the distinguished point for the 
entire cactus). 
As both maps respect degeneracies and face maps, 
it is straightforward to check that this induces a well-defined map 
\[ c : |S_{\bullet}|^n \rightarrow |S_{\bullet}| . \]
and to check that it gives an operad 
action on $|S_{\bullet}|$. 
\end{proof}

\vspace{6mm}

\section{The Hochschild cohomology of a Poincar\'{e} algebra}

Recall the following definition of Poincar\'{e} algebras, which is 
a Frobenius algebra that is also commutative.

\begin{defi}
An associative algebra $A$ is a Poincar\'{e} algebra if 
it is commutative, and has the property 
that there is an augmentation $A \rightarrow K$, where $K$ is the base field, 
such that the adjoint map 
\[ A \rightarrow A^{\vee} \]
obtained from
\[ A \otimes A \rightarrow A \rightarrow K \]
is an isomorphism. (Here, $A^{\vee} = Hom_{K}(A, K)$ is the 
dual of $A$ in the category of $K$-modules.)
\end{defi} 

In particular, by the isomorphism between 
$A$ and $A^{\vee}$, there is a coalgebra structure on $A$ which is 
dual to the algebra structure on $A$:
\begin{equation*}
 A \stackrel{\simeq}{\rightarrow} A^{\vee} \rightarrow 
(A \otimes A)^{\vee} \simeq A^{\vee} \otimes A^{\vee} 
\stackrel{\simeq}{\rightarrow} A \otimes A . 
\label{dualstr}
\end{equation*}
The motivating example is the cohomology of 
a orientable manifold $M$, which is isomorphic to its dual $H_{\ast}(M)$ by 
Poincar\'{e} duality. 

We consider the dual to $C_{cyclic}(A)$, which in this case is the
cyclic bar construction $B_{cyclic}(A)$  
(now considered as an algebra instead of a coalgebra). Note that 
for commutative $A$, this gives
the Hochschild homology $HH_{\ast}(A)$. The main theorem of this section is 
the following. 

\begin{theor}
Let $A$ be a Poincar\'{e} algebra. Then the dual of its Hochschild 
homology $HH_{\ast}(A)^{\vee}$ naturally has the structure of the 
an algebra over the chain complexes of $\mathcal{C}$. 
\end{theor}

By~(\ref{dualstr}), the algebra multiplication of $A$ 
is dual to the coproduct on $A$ as a coalgebra, and it is straightforward
to check that the algebra unit of $A$ is dual to the coalgebra counit of 
$A$. Hence, the degeneracies are given by the unit $K 
\rightarrow A$, and the face maps are given by multiplication $A \otimes A 
\rightarrow A$, and $B_{cyclic}(A)$ is the usual bar construction 
$B(A)$ of $A$, with 
an extra degeneracy. In fact, we can 
visualize $B^{n}_{cyclic}(A)$ as the tensor product of $n+1$ copies 
of $A$ arranged in a circle as below, with the appropriate 
coface and codegeneracy maps. (Note that it can also be considered as a 
cocyclic object, with the degeneracies being the cofaces, and the 
faces being the codegeneracies, since the categories of cyclic and cocyclic
objects are canonically equivalent.) 

\[ \begin{array}{ccccccc}
& & A & \otimes & A & & \\
& \otimes & & & & \otimes & \\
A & & & & & & A \\
\otimes & & & & & & \otimes \\
\vdots & & & & & & \vdots \\
& \otimes & & & & \otimes & \\
& & A & \otimes & A & &
\end{array} \]

We claim that this has the 
dual structure to a cactus object in the category of module over the 
base field $K$. (Note that here, the product $\times$ of sets is replaced by 
$\otimes$ of algebras). In other words, for every cactus configuration 
$(c, X)$, where there are $j_k$ points on the $k$-th lobe of the 
cactus, there is a map 
\begin{equation*}
B^{\Sigma j_k -1}_{cyclic}(A) \rightarrow B^{j_1 -1}_{cyclic}(A) 
\otimes \cdots \otimes B^{j_n -1}_{cyclic}(A) 
\label{cocactus}
\end{equation*}
and these maps commute with the face and degeneracy maps, as well as with
the composition of cacti. We will need the following lemma. 

\begin{lemma}
The cactus operad $\mathcal{C}$ is generated by 
$\mathcal{C}(1)$ and $\mathcal{C}(2)$.
\label{koszulgen}
\end{lemma}

\begin{proof}
For any cactus configuration $c \in \mathcal{C}(n)$, we can obtain 
$c$ from $S^1$ by a succession of operations that pinches one of the lobes
into two, which increases the number of lobes by one.
Note that each cactus $c$ determines a parametrization on $S^1$, 
which gives the 
parametrization on its lobes. Now $\mathcal{C}(1)$ consists of 
parametrized circles with one distinguished point. Hence, composition 
with elements of $\mathcal{C}(1)$ allows us to reparametrize $S^1$.
Thus, to obtain any $c \in \mathcal{C}(n)$, we can begin with the 
correct parametrization on
$S^1$. However, each pinching operation is precisely operad composition
with an element of $\mathcal{C}(2)$ in one position (and with the identity 
in $\mathcal{C}(1)$ in all other positions). 
Finally, any way of ordering the lobes 
can be obtained by the action of the symmetric group on $\mathcal{C}(n)$. 
\end{proof}

Given such a cactus configuration 
$(c, X)$, we consider $B_{cyclic}^{\Sigma j_i -1}(A)$. From 
the intersection points of $(c, X)$, 
we get a list of identifications in the set $\{ 1, \ldots, \Sigma j_k \}$. 
For each intersection point, we ``pinch'' $B_{cyclic}^{\Sigma j_i -1}$ 
by tensoring together the copies of $A$ at the positions to be 
identified, using the product $A \otimes A \rightarrow A$.
This gives a tensor 
product of copies of $A$, with 
one copy of $A$ corresponding to each point in $X$
in the configuration $c$, but with only one 
copy of $A$ at each intersection point. 
By Proposition~\ref{koszulgen}, it suffices to consider the 
situation when $c$ is a cactus 
with two lobes. In this case, the configuration will be of the form 
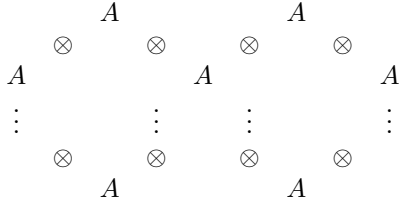
\[ \begin{array}{ccccccccc}
& & A & & & & A & & \\
& \otimes & & \otimes & & \otimes & & \otimes &\\
A & & & & A & & & & A \\
\vdots & & & \vdots & & \vdots & & & \vdots \\
& \otimes & & \otimes & & \otimes & & \otimes & \\
& & A & & & & A & &
\end{array} \]
To break this up into $A^{\otimes j_1 -1} \otimes \cdots \otimes 
A^{\otimes j_n-1}$, we use the coproduct structure 
$A \rightarrow A \otimes A$ to 
``break apart'' the circles: i.~e. we break apart the copy of $A$ 
at each intersection point into several copies, one for each lobe of
the cactus at that intersection point. Hence, the configuration will look 
like 
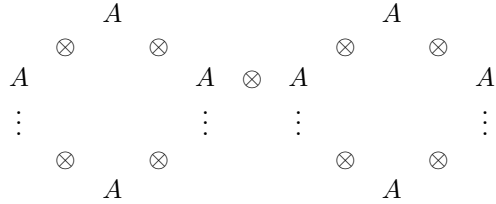
\[ \begin{array}{ccccccccccc}
& & A & & & & & & A & & \\
& \otimes & & \otimes & & & & \otimes & & \otimes & \\
A & & & & A & \otimes & A & & & & A \\
\vdots & & & & \vdots & & \vdots & & & & \vdots \\
& \otimes & & \otimes & & & & \otimes & & \otimes & \\
& & A & & & & & & A & & 
\end{array} \]

This gives maps of the form 
\[ B_{cyclic}^{j_1 + \cdots + j_n -1}(A) \rightarrow B_{cyclic}^{j_1 -1}(A)
\otimes \cdots \otimes B_{cyclic}^{j_n-1}(A)  . \]
Note that since $A$ is associative as an algebra
there is no ambiguity in the ordering when we are 
pinching together three or more copies of $A$ at an intersection point 
with three or more lobes. Further, since $A$ is also coassociative and 
cocommutative as a coalgebra, there is also no ambiguity
when breaking apart circles at an intersection point 
of three or more lobes of the cactus. Also, note that in order to consider 
each circle of $A$'s along a lobe of the cactus as $A^{\otimes j_k -1}$, 
we begin counting at the copy of $A$ at the position corresponding to the 
distinguished point of that lobe, and go in the 
direction of the parametrization of the 
lobe given in the cactus data. 

It remains to check that this is compatible with the 
faces and degeneracies, as well as the composition of operads. 
By Proposition~\ref{koszulgen},   
it suffices to consider the case of $\mathcal{C}(2)$, with only two copies
of $A$ being pinched together. It is clear that if 
the face or degeneracy map does not involve the $A$'s being pinched, there 
is no problem of compatibility with the faces or degeneracies. We need to 
consider the case when the face or degeneracy involves one of the copies 
of $A$ being pinched. For the degeneracies, we need to compare the following
maps, starting from $B_{cyclic}^n(A)$. First, we have the pinching followed by 
a degeneracy: 
\begin{equation*}
\begin{split}
{\begin{array}{ccccc}
& & A & & \\
& \otimes & & \otimes & \\
\vdots & & & & \vdots \\
& \otimes & & \otimes & \\
& & A & &
\end{array}} & \stackrel{\mu}{\rightarrow}
{\begin{array}{ccccccc}
& \cdots & \otimes & & \otimes & \cdots & \\
\vdots & & & A & & & \vdots \\
& \cdots & \otimes & & \otimes & \cdots & 
\end{array}} \stackrel{\psi}{\rightarrow} \\
{\begin{array}{ccccccc}
& \cdots & \otimes & & \otimes & \cdots & \\
\vdots & & & A \otimes A & & & \vdots \\
& \cdots & \otimes & & \otimes & \cdots & 
\end{array}} & \stackrel{\eta}{\rightarrow}
{\begin{array}{ccccccc}
& \cdots & \otimes & & \otimes A & \cdots & \\
\vdots & & & A \otimes A & & & \vdots \\
& \cdots & \otimes & & \otimes & \cdots & 
\end{array}} . 
\end{split}
\end{equation*}
Recal that here, $\mu$ is the product on $A$, and $\psi$ is the 
coproduct. These two together form the pinching operation. The 
last map $\eta$ is the unit $K \rightarrow A$, which inserts the 
new copy of $A$ in the circle on the right. 

On the other hand, taking the degeneracy first, and then pinching is 
\begin{equation*}
\begin{split}
{\begin{array}{ccccc}
& & A & & \\
& \otimes & & \otimes & \\
\vdots & & & & \vdots \\
& \otimes & & \otimes & \\
& & A & &
\end{array}} & \stackrel{\eta}{\rightarrow} 
{\begin{array}{ccccc}
& & A \otimes A & & \\
& \otimes & & \otimes & \\
\vdots & & & & \vdots \\
& \otimes & & \otimes & \\
& & A & & 
\end{array}} \stackrel{\mu}{\rightarrow} \\
{\begin{array}{ccccccc}
& \cdots & \otimes & & \otimes A & \cdots & \\
\vdots & & & A & & & \vdots \\
& \cdots & \otimes & & \otimes & \cdots & 
\end{array}} & \stackrel{\psi}{\rightarrow}
{\begin{array}{ccccccc}
& \cdots & \otimes & & \otimes A & \cdots & \\
\vdots & & & A \otimes A & & & \vdots \\
& \cdots & \otimes & & \otimes & \cdots & 
\end{array}} .
\end{split}
\end{equation*} 
The first map $\eta$ inserts the right one of the two $A$'s on the top 
of the second row. The product $\mu$ coming next multiplies together
the left copy of the two $A$'s on top with the copy of $A$ at the bottom of 
the circle, and the last map $\psi$ is the coproduct splitting the resulting 
single copy of $A$ into two. 
It is easy to see that these two compositions are equal to each other. 

For the face map, we need to compare the following maps, starting from 
$B_{cyclic}^n (A)$. To apply the face first, then pinch, we have
\begin{equation*}
\begin{split}
{\begin{array}{ccccc}
& & A \otimes A & & \\
& \otimes & & \otimes & \\
\vdots & & & & \vdots \\
& \otimes & & \otimes & \\
& & A & & 
\end{array}} & \stackrel{\mu}{\rightarrow}
{\begin{array}{ccccc}
& & A & & \\
& \otimes & & \otimes & \\
\vdots & & & & \vdots \\
& \otimes & & \otimes & \\
& & A & & 
\end{array}} \stackrel{\mu}{\rightarrow} \\
{\begin{array}{ccccccc}
& \cdots & \otimes & & \otimes & \cdots & \\
\vdots & & & A & & & \vdots \\
& \cdots & \otimes & & \otimes & \cdots & 
\end{array}} & \stackrel{\psi}{\rightarrow}
{\begin{array}{ccccccc}
& \cdots & \otimes & & \otimes & \cdots & \\
\vdots & & & A \otimes A & & & \vdots \\
& \cdots & \otimes & & \otimes & \cdots & 
\end{array}} .
\end{split}
\end{equation*}
The first $\mu$ is a face, multiplying together the two copies of $A$ at the top 
of the circle. The second $\mu$ and the coproduct $\psi$ together make up the 
pinching operation: the product $\mu$ multiplies together the copies of $A$ 
at the top and the bottom of the circle, and the $\psi$ splits the resulting 
single copy of $A$ into two copies. 

On the other hand, doing the pinching first, then the face gives
\begin{equation*}
\begin{split}
{\begin{array}{ccccc}
& & A \otimes A & & \\
& \otimes & & \otimes & \\
\vdots & & & & \vdots \\
& \otimes & & \otimes & \\
& & A & & 
\end{array}} & \stackrel{\mu}{\rightarrow}
{\begin{array}{ccccccc}
& \cdots & \otimes & & \otimes A & \cdots & \\
\vdots & & & A & & & \vdots \\
& \cdots & \otimes & & \otimes & \cdots & 
\end{array}} \stackrel{\psi}{\rightarrow} \\
{\begin{array}{ccccccc}
& \cdots & \otimes & & \otimes A & \cdots & \\
\vdots & & & A \otimes A & & & \vdots \\
& \cdots & \otimes & & \otimes & \cdots & 
\end{array}} & \stackrel{\mu}{\rightarrow} 
{\begin{array}{ccccccc}
& \cdots & \otimes & & \otimes & \cdots & \\
\vdots & & & A \otimes A & & & \vdots \\
& \cdots & \otimes & & \otimes & \cdots & 
\end{array}} .
\end{split}
\end{equation*}
Here, the first $\mu$ and the $\psi$ form the pinching operation, and the last 
$\mu$ is the face map. The first $\mu$ multiplies together the copies of $A$
at the top left and the bottom of the circle, and the coproduct $\psi$ splits
the resulting single $A$ (at the center of the configuration of two circles) 
into two copies. Finally, the last $\mu$ multiplies together one of 
the two resulting $A$'s at the center of the configuration (the one on the 
right) with the copy of $A$ immediately above it in the circle on the right. 
It is easy to see that by the associativity of $A$ as an algebra, the two 
compositions are equal. 

Thus,  
\[ B_{cyclic}(A) \]
has the dual structure~(\ref{cocactus}) to a cactus set. 
To turn it into a cactus set, we take its dual
\[ B_{cyclic}(A)^{\vee}= Hom_K(B_{cyclic}(A), K) . \]
This dual is taken termwise, so the $n$-th stage is $(A^{\vee})^{\otimes n+1}
\simeq A^{\otimes n+1}$, 
and the face and degeneracy maps are also dualized. 
As seen above, this is also the cyclic cobar construction 
$C_{cyclic}(A)$ on $A$ as a coalgebra, via the 
isomorphism $A \simeq A^{\vee}$, with a shift by the
dimension of $A$. This is still a cyclic 
module (since the dual of a cyclic module is a cyclic module), and now the cactus 
structure maps go in the right direction:
\[ C_{cyclic, j_1 -1}(A) \otimes \cdots \otimes 
C_{cyclic, j_n -1}(A) \rightarrow C_{cyclic, j_1 + \cdots + j_n -1}(A). \]

Hence, applying the totalization functor, by the result from 
the previous section, 
we get an action of the chain complexes cactus operad 
$\mathcal{C}$ on the dual of the Hochschild homology
$HH_{\ast}(A)^{\vee}$. Voronov's theorem tells us that the 
cactus operad is equivalent to $\mathcal{D}_2$. 

To get to the Hochschild cohomology $HC^{\ast}(A)$ of $A$, recall that
for a commutative algebra $A$,
$HC^{\ast}(A)$ can be calculated as the cohomology of 
the cyclic 
cobar construction $C_{cyclic}(A)$ of $A$~(\cite{hu1}, Proposition 2.8). (This is 
the same as the usual cobar construction $C(A)= C_{K}(K, A, K)$ over $K$
of $A$ as a coalgebra, 
i.~e. $C^{n}_{cyclic}(A) = A^{\otimes n+1}$, with 
coface maps given by the counit on $A$ as 
a coalgebra (which is the same as 
the augmentation $A \rightarrow K$), and codegeneracy maps given 
by multiplying together two adjacent copies of $A$. The 
only difference is that there is an extra codegeneracy between 
the first and last copies of $A$ at each stage, and that there is also 
a rotational action on the $n$-th stage by $\mathbb{Z}/({n+1})$.)
Thus, for a Poincar\'{e} algebra 
$A$, the dual of the Hochschild homology gives $HC^{\ast}(A)$.
\vspace{3mm}

\noindent\textbf{Remark:}
The motivating example for $A$ is the cochain complex $C^{\ast}(M)$ for 
$M$ a smooth compact orientable manifold, given in Chas and Sullivan~\cite{cs}. 
However, we must note that this case is only an example in the philosophical
sense, since $C^{\ast}(M)$ is only a Poincar\'{e} algebra up to homotopy, 
which we do not address in this paper. By Theorem 3 of~\cite{cj}, 
\[ HC^{\ast}(C^{\ast}(M)) \simeq C_{\ast}(LM^{\nu_M}) \]
where $\nu_M$ is the stable normal bundle of $M$. This corresponds to the 
$\mathcal{D}_2$-algebra structure on $H_{\ast}(LM)$ given in~\cite{cs}, 
and the shift is by the dimension of the manifold $M$.


\begin{thebibliography}{99}


\bibitem{bhm} M. Bokstedt, W. C. Hsiang and I. Madsen. The cyclotomic trace
and algebraic $K$-theory of spaces. \emph{Invent. Math.} 111 (1993), no. 3, 
465-539. 

\bibitem{connes} A. Connes. Cohomologie cyclique et Foncteurs 
$Ext^n$. \emph{C. R. Acad. Sci. Paris} 296 (1983) 953-958. 

\bibitem{cs} M. Chas and D. Sullivan. String topology. math.GT/9911159.

\bibitem{cj} R. Cohen and J. D. S. Jones. A homotopy theoretic realization of 
string topology. math.GT/0107187.

\bibitem{elmendorf} A. D. Elmendorf. A simple formula for cyclic duality. 
\emph{Proc. Amer. Math. Soc.} 118 (1993) 3, 709-711. 

\bibitem{hu1} P. Hu. Higher string topology on general spaces. To appear in 
\emph{Proc. London Math. Soc.}

\bibitem{jones} J. D. S. Jones. Cyclic homology and equivariant homology. 
\emph{Invent. Math.} 87 (1987) no. 2, 403-423.


\bibitem{ms} J. E. McClure and J. Smith. A solution of Deligne's hochschild 
cohomology conjecture. \emph{Contemp. Math.} 293, 153-193. Amer. Math. Soc., 2002.

\bibitem{voronov} A. A. Voronov. Notes on universal algebras. QA/0111009. 

\end{thebibliography}
\end{document}